\documentclass[a4paper, 14pt]{extarticle}
\usepackage{amsfonts,amsmath,amssymb,amsthm,amscd}
\usepackage[T2A]{fontenc}
\usepackage[cp1251]{inputenc}
\usepackage[english]{babel}
\usepackage{graphicx}
\usepackage[papersize={210mm, 297mm},left=2cm, right=2cm, top=2cm, bottom=2cm]{geometry}

\numberwithin{equation}{section}

\theoremstyle{definition}
\newtheorem{definition}{Definition}[section]

\newtheorem{example}[definition]{Example}

\theoremstyle{theorem}
\newtheorem{theorem}[definition]{Theorem}

\newtheorem{corollary}[definition]{Corollary}

\numberwithin{equation}{section}

\begin{document}

\noindent UDK 515.124

\smallskip

\noindent MATHEMATICS

\smallskip

\noindent Institute of Applied Mathematics and Mechanics NAS of Ukraine, Dobrovolskogo~Str.~1, Sloviansk, 84100, Ukraine; Phone: +38(0626)66 55 00

\smallskip

\noindent \emph{Emails:} viktoriiabilet@gmail.com, oleksiy.dovgoshey@gmail.com

\begin{center}
\large\textbf{Asymptotic behavior of metric spaces at infinity}
\end{center}

\begin{center}
\textbf{Viktoriia Bilet, Oleksiy Dovgoshey}

\end{center}

\smallskip


\medskip

\noindent\emph{A new sequential approach to investigations of structure of metric spaces at infinity is proposed. Criteria for finiteness and boundedness of metric spaces at infinity are found.}

\medskip


\section{Introduction}

Under the asymptotic structure describing the behavior of an unbounded metric space $(X,d)$ at infinity we mean a metric space that is a limit of rescaling metric spaces $(X, \frac{1}{r_n} d)$ for $r_n$ tending to infinity. The Gromov--Hausdorff convergence and the asymptotic cones are most often used for construction of such limits. Both of these approaches are based on higher-order logic abstractions (see, for example, \cite{Ro} for details), which makes them very powerful, but it does away the constructiveness. In this paper we propose a more elementary, sequential approach for describing the structure of unbounded metric spaces at infinity.

Let $(X,d)$ be an unbounded metric space, $p$ be a point of $X$ and $\tilde r = (r_n)_{n\in\mathbb N}$ be a scaling sequence of positive real numbers tending to infinity. Denote by $\mathbf{\tilde X}_{\infty, \tilde r}$ the set of all sequences $\tilde x=(x_n)_{n\in\mathbb N}\subset X$ for each of which $\mathop{\lim}\limits_{n\to\infty}d(x_n, p)=\infty$ and there is a finite limit $\tilde{\tilde d}_{\tilde r}(\tilde x) := \mathop{\lim}\limits_{n\to\infty}\frac{d(x_n, p)}{r_n}$. Define the equivalence relation $\equiv$ on $\mathbf{\tilde X}_{\infty, \tilde r}$ as
\begin{equation*}
(\tilde x\equiv\tilde y)\Leftrightarrow\left(\lim_{n\to\infty}\frac{d(x_n, y_n)}{r_n}=0\right).
\end{equation*}
Let $\mathbf{\Omega}_{\infty, \tilde r}^{X}$ be the set of equivalence classes generated by $\equiv$ on $\mathbf{\tilde X}_{\infty, \tilde r}$. We shall say that points $\alpha, \beta\in\mathbf{\Omega}_{\infty, \tilde r}^{X}$ are \emph{mutually stable} if for $\tilde x\in\alpha$ and $\tilde y\in\beta$ there is a limit
\begin{equation}\label{e1.1}
\rho(\alpha, \beta):=\lim_{n\to\infty}\frac{d(x_n, y_n)}{r_n}.
\end{equation}
Let us consider the weighted graph $(G_{X, \tilde r}, \rho)$ with the vertex set $V(G_{X, \tilde r})=\mathbf{\Omega}_{\infty, \tilde r}^{X}$, and the edge set $E(G_{X, \tilde r})$ such that
$$
\left(\{u, v\}\in E(G_{X, \tilde r})\right)\Leftrightarrow\left(u \,\mbox{and}\, v \, \mbox{are mutually stable and}\, u\ne v\right),
$$
and the weight $\rho:E(G_{X, \tilde r})\to\mathbb R^{+}$ defined by formula~\eqref{e1.1}.

\begin{definition}
The \emph{pretangent spaces} (to $(X, d)$ at infinity w.r.t. $\tilde r$) are the maximal cliques $\Omega_{\infty,\tilde{r}}^X$ of $G_{X, \tilde r}$ with metrics determined with the help of~\eqref{e1.1}.
\end{definition}

Recall that a \emph{clique} in a graph $G$ is a set $A\subseteq V(G)$ such that every two distinct points of $A$ are adjacent. A clique $C$ in $G$ is maximal if $C\subseteq A$ implies $C=A$ for every clique $A$ in $G$.

Define the subset $\alpha_0 = \alpha_0(X, \tilde{r})$ of the set $\tilde{\mathbf{X}}_{\infty, \tilde{r}}$ by the rule:
\begin{equation}\label{e1.2}
\left(\tilde{z} \in \alpha_0\right) \Leftrightarrow \left(\tilde{z} \in \tilde{\mathbf{X}}_{\infty, \tilde{r}} \text{ and } \tilde{\tilde{d}}_{\tilde{r}} (\tilde{z}) = 0\right),
\end{equation}
Then $\alpha_0$ is a common point of all pretangent spaces $\Omega_{\infty,\tilde{r}}^X$. It means in particular that the graph $G_{X, \tilde{r}}$ is connected.

Let $(n_k)_{k\in\mathbb N}\subset\mathbb N$ be infinite and strictly increasing. Denote by $\tilde r'$ the subsequence $(r_{n_k})_{k\in \mathbb N}$ of the $\tilde r=(r_n)_{n\in\mathbb N}$ and, for every $\tilde x=(x_n)_{n\in\mathbb N}\in\mathbf{\tilde X}_{\infty, \tilde r}$, write $\tilde x' := (x_{n_k})_{k\in\mathbb N}$. It is clear that $\mathop{\lim}\limits_{k \to \infty} d(x_{n_k},p) = \infty$ and $\tilde{\tilde d}_{\tilde r'}(\tilde x') = \tilde{\tilde d}_{\tilde r}(\tilde x)$ for every $\tilde x\in\mathbf{\tilde X}_{\infty, \tilde r}$. Moreover, if $\tilde y\in\mathbf{\tilde X}_{\infty, \tilde r}$ and $\mathop{\lim}\limits_{n\to\infty}\frac{d(x_n, y_n)}{r_n}$ exists, then
\begin{equation}\label{e1.3}
\lim_{k\to\infty}\frac{d(x_{n_k}, y_{n_k})}{r_{n_k}}=\lim_{n\to\infty}\frac{d(x_n, y_n)}{r_n}.
\end{equation}
Let us define $\tilde{\mathbf{X}}_{\infty, \tilde{r}'}$ and $\rho'$ similarly to $\tilde{\mathbf{X}}_{\infty, \tilde{r}}$ and, respectively, $\rho$ and let
$\pi\colon \tilde{\mathbf{X}}_{\infty, \tilde{r}} \to \mathbf{\Omega}_{\infty, \tilde{r}}$, $\pi'\colon \tilde{\mathbf{X}}_{\infty, \tilde{r}} \to \mathbf{\Omega}_{\infty, \tilde{r}'}$ be the natural projections
$\pi(x): = \{\tilde{y} \in \tilde{\mathbf{X}}_{\infty, \tilde{r}}\colon \rho(\tilde{x}, \tilde{y}) = 0\}$, $\pi'(x): = \{\tilde{y} \in \tilde{\mathbf{X}}_{\infty, \tilde{r}'}\colon \rho'(\tilde{x}, \tilde{y}) = 0\}$ and let $\varphi_{\tilde{r}'}(\tilde{x}):= \tilde{x}'$ for all $\tilde{x} \in \tilde{\mathbf{X}}_{\infty, \tilde{r}}$. Then there is an embedding $em' \colon \mathbf{\Omega}_{\infty, \tilde{r}}^X \to \mathbf{\Omega}_{\infty, \tilde{r}}^X$ of the weighted graph $(G_{X, \tilde{r}}, \rho)$ in the weighted graph $(G_{X, \tilde{r}'}, \rho')$ such that the diagram
$$
\begin{CD}
\tilde{\mathbf{X}}_{\infty, \tilde r} @>\varphi_{\tilde r'}>> \tilde{\mathbf{X}}_{\infty, \tilde r'}\\
@V{\pi}VV @VV{\pi'}V \\
\mathbf{\Omega}_{\infty, \tilde r}^{X} @>em'>> \mathbf{\Omega}_{\infty, \tilde r'}^{X}
\end{CD}
$$
is commutative. Since $em'$ is an embedding of weighted graphs, $em'(C)$ is a clique in $G_{X, \tilde{r}'}$ if $C$ is a clique in $G_{X, \tilde{r}}$. Furthermore, \eqref{e1.3} implies that the restrictions $em'|_{\Omega_{\infty, \tilde r}^{X}}$ are isometries of the pretangent spaces $\Omega_{\infty, \tilde r}^{X}$ on the metric spaces $em'(\Omega_{\infty, \tilde r}^{X})$.

\begin{definition}
A pretangent space $\Omega_{\infty, \tilde r}^{X}$ is \emph{tangent} if the clique $em'(\Omega_{\infty, \tilde r}^{X})$ is maximal for every infinite, strictly increasing sequence $(n_k)_{k\in\mathbb N}\subset\mathbb N$.
\end{definition}

\begin{example}
Let $E$ be a finite-dimensional Euclidean space and let $X \subseteq E$ such that the Hausdorff distance $d_H(E, X)$ is finite. Then for every scaling sequence $\tilde{r}$ all pretangent spaces $\Omega_{\infty, \tilde r}^X$ are tangent and isometric to $E$.
\end{example}

In conclusion of this brief introduction it should be noted that there exist other techniques which allow to investigate the asymptotic properties of metric spaces at infinity. As examples, we mention only the balleans theory~\cite{PZ} and the Wijsman convergence \cite{LechLev}, \cite{Wijs64}, \cite{Wijs66}.

\section{Finiteness}

In this section we study the conditions under which pretangent spaces are finite.

\begin{theorem}\label{th2.1}
Let $(X, d)$ be an unbounded metric space, $p\in X$, $n\ge 2$ and let
\begin{equation*}\label{Fn}
F_n(x_1, ..., x_n):=
\begin{cases}
\frac{\mathop{\min}\limits_{1\le k\le n}d(x_k, p)\mathop{\prod}\limits_{1\le k<l\le n}d(x_k, x_l)} {\left(\mathop{\max}\limits_{1\le k\le n}d(x_k, p)\right)^{\frac{n(n-1)}{2}+1}}, & \mbox{if} $ $ (x_1, ..., x_n) \ne (p, ..., p)\\
0, & \mbox{if}$ $ (x_1, ..., x_n) = (p, ..., p). \\
\end{cases}
\end{equation*}
Then the inequality $\left|\Omega_{\infty, \tilde r}^{X}\right|\le n$ holds for every pretangent space $\Omega_{\infty, \tilde r}^{X}$ if and only if $\mathop{\lim}\limits_{x_1, ..., x_n\to\infty}F_{n}(x_1, ..., x_n)=0.$
\end{theorem}

Note that, for every unbounded metric space $(X, d)$, there is a pretangent space $\Omega_{\infty, \tilde r}^{X}$ consisting at least two points.

\begin{corollary}
Let $(X, d)$ be an unbounded metric space and let $\alpha_0$ be a point defined by \eqref{e1.2} for every scaling sequence $\tilde r.$ Then the following statements are equivalent.
\begin{enumerate}
\item[\rm(i)] The graph $G_{X, \tilde r}$ is a star with the center $\alpha_0$ for every scaling sequence $\tilde r$;
\item[\rm(ii)] The limit relation $\mathop{\lim}\limits_{x_1, x_2\to\infty}F_{2}(x_1, x_2)=0$ holds.
\end{enumerate}
\end{corollary}

Let us consider now the problem of existence of finite \emph{tangent} spaces.

\begin{definition}\label{d2.2}
Let $E\subseteq\mathbb R^{+}$. \emph{The porosity} of $E$ at infinity is the quantity
\begin{equation}\label{infpor}
p(E, \infty):=\limsup_{h\to\infty}\frac{l(\infty, h, E)}{h}
\end{equation}
where $l(\infty, h, E)$ is the length of the longest interval in the set $[0, h]\setminus E.$ The set $E$ is \emph{strongly porous} at infinity if $p(E, \infty)=1$.
\end{definition}

The standard definition of the porosity at a point can be found in~\cite{Th}.

For a metric space $(X, d)$ and $p\in X$ write $S_{p}(X):=\{d(x, p)\colon x\in X\}$.

\smallskip

\begin{theorem}\label{th2.3}
Let $(X, d)$ be an unbounded metric space, $p\in X.$ The following statements are equivalent:
\begin{enumerate}
\item[\rm(i)] The set $S_{p}(X)$ is strongly porous at infinity;
\item[\rm(ii)] There is a single-point tangent space $\Omega_{\infty, \tilde r}^{X}$;
\item[\rm(iii)] There is a finite tangent space $\Omega_{\infty, \tilde r}^{X}$;
\item[\rm(iv)] There is a compact tangent space $\Omega_{\infty, \tilde r}^{X}$;
\item[\rm(v)] There is a bounded, separable tangent space $\Omega_{\infty, \tilde r}^{X}$.
\end{enumerate}
\end{theorem}

Some results which are similar to Theorem~\ref{th2.1} and Theorem~\ref{th2.3} can be found in~\cite{ADK} and~\cite{DM} respectively.

\section{Boundedness}

Let $\tilde \tau=(\tau_n)_{n\in\mathbb N} \subset \mathbb{R}$. We shall say that $\tilde \tau$ is eventually increasing if the inequality $\tau_{n+1}\ge\tau_{n}$ holds for sufficiently large $n$. For $E\subseteq\mathbb R^{+}$ write $\tilde E_{\infty}^{i}$ for the set of eventually increasing sequences $\tilde \tau \subset E$ with $\mathop{\lim}\limits_{n\to\infty}\tau_{n}=\infty$. Denote also by $\tilde I_{E}^{i}$ the set of all sequences of open intervals $(a_n,b_n)\subseteq \mathbb R^{+}$ meeting the following conditions:

$\bullet$ Each $(a_n, b_n)$ is a connected component of the set $\mathop{Int}(\mathbb R^{+} \setminus E);$

$\bullet$ $(a_n)_{n\in\mathbb N}$ is eventually increasing;

$\bullet$ $\mathop{\lim}\limits_{n\to\infty}a_{n}=\infty$ and $\mathop{\lim}\limits_{n\to\infty}\frac{b_n-a_n}{b_n}=1$.

\medskip

Define an equivalence $\asymp$ on the set of sequences of strictly positive numbers as follows. Let $\tilde a=(a_n)_{n\in\mathbb N}$ and $\tilde{\gamma}=(\gamma_n)_{n\in\mathbb N}$. Then $\tilde a \asymp \tilde {\gamma}$ if there are some constants $c_1, c_2 >0$ such that $c_1 a_n < \gamma_n < c_2 a_n$ for every $n\in\mathbb N$.

\begin{definition}
Let $E\subseteq\mathbb R^{+}$ and let $\tilde \tau \in \tilde E_{\infty}^{i}$. The set $E$ is $\tilde \tau$-\emph{strongly porous} at infinity if there is a sequence $((a_n, b_n))_{n\in\mathbb N}\in\tilde I_{E}^{i}$ such that $\tilde\tau \asymp \tilde a$ where $\tilde a=(a_n)_{n\in\mathbb N}$. The set $E$ is \emph{completely strongly porous} at infinity if $E$ is $\tilde \tau$-strongly porous at infinity for every $\tilde \tau \in \tilde E_{\infty}^{i}$.
\end{definition}

Note that every completely strongly porous at infinity set is strongly porous at infinity but not conversely.

\begin{definition}
Let $(X,d)$ be an unbounded metric space and let $p\in X$. A scaling sequence $\tilde{r}$ is \emph{normal} if $\tilde{r}$ is eventually increasing and there is $\tilde{x} \in \tilde{\mathbf{X}}_{\infty, \tilde{r}}$ such that
\begin{equation*}\label{L4.4*}
\lim_{n\to\infty} \frac{d(x_n,p)}{r_{n}}=1.
\end{equation*}
\end{definition}

Write $\mathcal{F}_n(X)$ for the set of all pretangent spaces $\Omega_{\infty,\tilde r}^{X}$ with normal scaling sequences $\tilde{r}$. Under what conditions the family $\mathcal{F}_n(X)$ is uniformly bounded?

Recall that a family $\mathcal{F}$ of a metric spaces $(Y, d_Y)$ is \emph{uniformly bounded} if
$$
\sup_{Y \in \mathcal{F}} \mathop{\mathrm{diam}} Y < \infty.
$$
If all metric spaces $(Y, d_Y)$ are pointed with marked points $p_Y \in Y$ and
$$
\inf_{Y \in \mathcal{F}} \inf\{d_Y(p_Y, y)\colon y \in Y\setminus \{p_Y\}\} > 0,
$$
then we say that $\mathcal{F}$ is \emph{uniformly discrete} (w.r.t. the points $p_Y$).

The following theorem is an analog of Theorem~3.11 and Theorem~4.1 from~\cite{BDK}.

\begin{theorem}
Let $(X,d)$ be an unbounded metric space and let $p\in X$. Then the following statements are equivalent
\begin{enumerate}
\item [$(i)$] The family $\mathcal{F}_n(X)$ is uniformly bounded.
\item [$(ii)$] $S_{p}(X)$ is completely strongly porous at infinity.
\item [$(iii)$] The family $\mathcal{F}_n(X)$ is uniformly discrete w.r.t. the points $\alpha_0$ defined by~\eqref{e1.2}.
\end{enumerate}
\end{theorem}

If $\mathcal{F}_n(X)$ is uniformly bounded, then every pretangent space $\Omega_{\infty, \tilde r}^X$ is bounded, but the converse, in general, does not hold.

\begin{definition}
The set $E \subseteq \mathbb{R}^+$ is $\omega$-strongly porous at infinity if for every sequence $\tilde{\tau}\in \tilde{E}_{\infty}^i$ there is a subsequence $\tilde{\tau}'$ for which $E$ is $\tilde{\tau}'$-strongly porous at infinity.
\end{definition}

The following theorem gives a boundedness criterion for pretangent spaces.
\begin{theorem}\label{unbpret}
Let $(X,d)$ be an unbounded metric space and let $p \in X$. All pretangent spaces to $X$ at infinity are bounded if and only if the set $S_{p}(X)$ is $\omega$-strongly porous at infinity.
\end{theorem}

The example of $\omega$-strongly porous at infinity set $E \subseteq \mathbb{R}^+$ which is not completely strongly porous at infinity can be obtained as a modification of Example~2.10~\cite{BD}.

Using Theorem~\ref{unbpret} we can obtain a criterion of existence of an unbounded pretangent space. The condition of such type are importnt for development of a theory of pretangent spaces of the second order or more (i.e., pretangent spaces to pretangent spaces, pretangent spaces to pretangent spaces to pretangent spaces and so on).

Let $(Y, d_Y)$ be a metric space. Then for every $K\subseteq Y$ and $y\in Y$ we write $\textrm{dist}(y, K)=\inf\{d_{Y}(y, x): x\in K\}.$ The following definition can be found in \cite{AF}.

\begin{definition}
Let $(K_n)_{n\in\mathbb N}$ be a sequence of subsets of $(Y, d_Y).$ The set
$$\textrm{Lim}\inf_{n\to\infty}K_{n}:=\{y\in Y: \lim_{n\to\infty}\textrm{dist}(y, K_n)=0\}$$
is the \emph{Kuratowski lower limit} of $(K_n)_{n\in\mathbb N}$ in $(Y, d_Y).$
\end{definition}

For $A\subseteq\mathbb R$ and $t\in\mathbb R$ we set $tA:=\{ta: a\in A\}.$

\begin{theorem}\label{exunb}
Let $(X, d)$ be an unbounded metric space and let $p\in X.$ Then the following statements are equivalent.
\begin{enumerate}
\item [$(i)$] There exists an unbounded pretangent space $\Omega_{\infty, \tilde r}^{X}.$
\item [$(ii)$] There exists a scaling sequence $\tilde r=(r_n)_{n\in\mathbb N}$ such that the Kuratowski lower limit $\mathop{\emph{Lim}\inf}\limits_{n\to\infty}\frac{1}{r_n}S_{p}(X)$ is an unbounded subset of $\mathbb R^{+}.$
\item [$(iii)$] The set $S_{p}(X)$ is not $w$-strongly porous at infinity.
\end{enumerate}
\end{theorem}

Let $(\tilde r_m)_{m\in\mathbb N}$ be a sequence of scaling sequences. Then for every $m\in\mathbb N$ and every unbounded metric space $(X, d)$ we define a pretangent space $\Omega_{\infty, (\tilde r_1, ..., \tilde r_m)}^{X}$ by the following inductive rule: $\Omega_{\infty, (\tilde r_1)}^{X}:=\Omega_{\infty, \tilde r_1}^{X}$ if $m=1$ and $$\Omega_{\infty, (\tilde r_1, ..., \tilde r_m)}^{X}:=\Omega_{\infty, \tilde r_m}^{\Omega_{\infty, (r_1, ..., r_{m-1})}^{X}}$$ if $m\ge 2.$

Using Theorem~\ref{exunb} we can obtain the following corollary.

\begin{corollary}
Let $(X, d)$ be an unbounded metric space and let $p\in X.$ If the equality $p(S_{p}(X), \infty)=0$ holds, then there is a sequence $(\tilde r_m)_{m\in\mathbb N}$ of scaling sequences such that the pretangent space $\Omega_{\infty, (\tilde r_1, ..., \tilde r_m)}^{X}$ is unbounded for every $m\in\mathbb N.$
\end{corollary}

\noindent\textbf{Acknowledgments.} The research was supported by the grant of the State Fund for Fundamental Research (project F71/20570).

\small

\noindent Інститут прикладної математики і мехініки НАН України, вул. Добровольського, 1, 84100, м. Слов'янськ, Україна; Тел: +38(0626)66 55 00

\noindent E-mails: viktoriiabilet@gmail.com, oleksiy.dovgoshey@gmail.com

\medskip

\noindent \textbf{В.~В. Білет, О.~А. Довгоший}

\medskip

\noindent \textbf{Асимптотична поведінка метричних просторів на нескінченності}

\medskip

\noindent \emph{Запропоновано новий секвенціальний підхід до дослідження структури метричних просторів у нескінченно віддаленій точці. Знайдено критерії скінченності та обмеженості метричних просторів на нескінченності.}

\medskip

\noindent \textbf{Ключові слова:} асимптотична обмеженість метричного простору, асимптотична скінченність метричного простору, збіжність метричних просторів, сильна пористість у точці.

\bigskip

\noindent Институт прикладной математики и механики НАН Украины, ул. Добровольского, 1, 84100, г. Славянск, Украина; Тел: +38(0626)66 55 00

\noindent E-mails: viktoriiabilet@gmail.com, oleksiy.dovgoshey@gmail.com

\medskip

\noindent \textbf{В.~В. Билет, А.~А. Довгошей}

\medskip

\noindent \textbf{Асимптотическое поведение метрических пространств на бесконечности}

\medskip

\noindent \emph{Предложен новый секвенциальный подход к исследованию структуры метрических пространств в бесконечно удаленной точке. Найдены критерии конечности и ограниченности метрических пространств на бесконечности.}

\medskip

\noindent \textbf{Ключевые слова:} асимптотическая ограниченность метрического пространства, асимптотическая конечность метрического пространства, сходимость метрических пространств, сильная пористость в точке.


\begin{thebibliography}{10}
\bibitem{Ro} \emph{Roe J.} Lectures on coarse geometry. -- University Lecture Series 31, Providence, RI: American Mathematical Society, 2003.

\bibitem{PZ} \emph{Protasov I., Zarichnyi M.} General Asymptology. -- Mathematical Studies Monograph Series, \textbf{12}. L'viv: VNTL Publishers, 2007.

\bibitem{LechLev} \emph{Lechiki~A., Levi~S.} Wijsmann convergence in the hyperspace of a metric space // Bull. Unione Mat. Ital. -- 1987. -- 1-B~. -- P. 439--452.

\bibitem{Wijs64} \emph{Wijsman R.~A.} Convergence of sequences of convex sets, cones and functions // Bull. Amer. Math. Soc. -- 1964. -- \textbf{70}. -- P. 186--188.

\bibitem{Wijs66} \emph{Wijsman R.~A.} Convergence of sequences of convex sets, cones and functions II. // Trans. Amer. Math. Soc. -- 1966. -- \textbf{123}, No 1. -- P.~32--45.

\bibitem{Th} \emph{Thomson B.~S.} Real Functions. -- Lecture Notes in Mathematics, 1170. -- Springer-Verlag, Berlin, Heidelberg, New York, Tokyo, 1985.

\bibitem{ADK} \emph{Abdullayev F., Dovgoshey O., K\"{u}\c{c}\"{u}kaslan M.} Compactness and boundedness of tangent spaces to metric spaces // Beitr. Algebra Geom. -- 2010. -- \textbf{51}, No 2. -- P.547--576.

\bibitem{DM} \emph{Dovgoshey O., Martio O.} Tangent spaces to general metric spaces // Rev. Roumaine Math. Pures. Appl. -- 2011. -- \textbf{56}, No 2. -- P. 137--155.

\bibitem{BDK} \emph{Bilet V., Dovgoshey O.,  K\"{u}\c{c}\"{u}kaslan M.} Uniform boundedness of pretangent spaces, local constancy of metric derivatives and strong right upper porosity at a point // J. Analysis. -- 2013. -- \textbf{21}. -- P. 31--55.

\bibitem{BD} \emph{Bilet V., Dovgoshey O.} Boundedness of pretangent spaces to general metric spaces // Ann. Acad. Sci. Fenn. -- 2014. -- \textbf{39}, No 4. -- P.~73--82.

\bibitem{AF} \emph{Aubin J.-P., Frankowska H.} Set-Valued Analysis. -- Birkh\"{a}user, Boston, Basel, Berlin, 1990.


\end{thebibliography}
\end{document}